\documentclass[12pt,twoside]{article} 
\usepackage{amsfonts}
\usepackage{graphicx} 

\newtheorem{theorem}{Theorem}[section]
\newtheorem{lemma}[theorem]{Lemma}
\newtheorem{proposition}[theorem]{Proposition}

\newtheorem{definition}[theorem]{Definition}
\newtheorem{rmrk}[theorem]{Remark}

\newcommand{\R} {{\mathbb R}}

\newcommand{\N} {{\mathbb N}}
\newcommand{\qed} {\hfill {\small Q.E.D.} \par\medskip}
\newcommand{\skippar} {\par\medskip}
\newcommand{\ds} {\displaystyle}
\newcommand{\proof} {\noindent \textsc{Proof.} }
\newcommand{\proofof}[1] {\noindent \textsc{Proof of {#1}.} }
\newcommand{\article}[3] {\textsc{{#1}}, {\itshape {#2}}, {{#3}}.}
\newcommand{\book}[3] {\textsc{{#1}}, {\itshape {#2}}, {{#3}}.}
\newcommand{\vol} {\textbf}
\newcommand{\eps} {\varepsilon}
\newcommand{\rset}[2] {\left\{ #1 \: \left| \: #2 \right. \! \right\} }

\newcommand{\fn} {function}
\newcommand{\bi} {billiard}
\newcommand{\me} {measure}
\newcommand{\tr} {trajector}
\newcommand{\erg} {ergodic}
\newcommand{\sy} {system}
\newcommand{\hyp} {hyperbolic}
\newcommand{\sca} {scatterer}
\renewcommand{\o} {orbit}

\newcommand{\si} {\mathcal{S}}  
\newcommand{\Sc} {\mathcal{O}}  
\renewcommand{\a} {\alpha}      
\renewcommand{\b} {\beta}       
\newcommand{\is} {\mathcal{I}}  
\newcommand{\ij} {\mathcal{J}}  
\newcommand{\nh} {\mathcal{N}}  
\newcommand{\ps} {\mathcal{M}}  
\newcommand{\x} {x}             
\newcommand{\y} {y}             
\newcommand{\ts} {\mathrm{T}}   
\newcommand{\co} {\mathcal{C}}  
\newcommand{\ph} {\varphi}      
\newcommand{\wu} {W^u}          
\newcommand{\ws} {W^s}          
\newcommand{\wsu} {W^{s(u)}}    

\newcommand{\sect}[1] {\section{{#1}} \setcounter{equation}{0}}

\newcommand{\fig}[3] {
\medskip\smallskip
\begin{figure}[ht]
        \centering
        \includegraphics[width=#2]{#1.eps}
        \begin{minipage}[t]{0.75\linewidth} 
                \caption{\baselineskip=14pt {#3}}
                \protect\label{#1}
        \end{minipage}
\end{figure}
\medskip
}

\newenvironment{remark}
{\begin{rmrk} \em}
{\end{rmrk}}

\begin{document}

\title{Aperiodic Lorentz gas: recurrence and ergodicity}

\author{
        \textsc{Marco Lenci} \\
        Department of Mathematical Sciences \\
        Stevens Institute of Technology \\
        Hoboken, NJ \ 07030, \ U.S.A. \\
        \footnotesize{E-mail: \texttt{mlenci@stevens-tech.edu}}
}

\date{June 2002}

\maketitle

\begin{abstract}
        We prove that any generic (i.e., possibly aperiodic) Lorenz
        gas in two dimensions, with finite horizon and non-degenerate
        geometrical features, is ergodic if it is recurrent. We also
        give examples of aperiodic recurrent gases.

\end{abstract}

\sect{Introduction}
\label{sec-intro}

By Lorentz gas (LG) we mean the free motion of a point particle in the
plane subject to elastic collisions against a fixed array of
dispersing \sca s. Each \sca\ $\Sc_\a$ is an open, bounded, connected,
simply connected, strictly convex domain of $\R^2$, with smooth
boundary, and is labeled by the variable $\a \in \is$. The \sca s are
assumed to be pairwise disjoint. A LG is called \emph{periodic} when
$\cup_{\a \in \is} \Sc_\a$ is left invariant by $G$, a discrete group
of translations in the plane, whose fundamental domain $\R^2 / G$ is
compact.

Notice the departure from the usual terminology, whereby `Lorentz gas'
means `periodic Lorentz gas' (PLG). To avoid possible confusions let
us also recall that some literature calls `Lorentz gas' the dynamical
\sy\ defined on $\R^2 / G$, in the periodic case. However, common
practice has nowadays christened such \sy s \emph{Sinai \bi s}.

\skippar

The LG has a long and honorable history. It was first introduced in
its periodic and three-dimensional version by Lorentz in 1905, to
describe the motion of an electron in a crystal \cite{lo}.  Later in
the century it became extremely popular in statistical mechanics and
\erg\ theory. The scientific community was seeking a rather realistic
model whose stochastic properties could be rigorously proven---in
particular Boltzmann's Ergodic Hypothesis (namely, \erg ity). It is
because of this connection with statistical mechanics that it was
renamed `Lorentz \emph{gas}'. At this point the periodicity condition
had lost its \emph{raison d'\^etre}, but still remained in virtually
all the publications (at least as far as I know), on grounds of
mathematical convenience.

As a matter of fact, most of the massive literature in this area of
research is actually concerned with Sinai \bi s, and many results on
the PLG were generated as corollaries of theorems on the Sinai \bi.
It is definitely outside the scope of this introduction to mention all
that is known about these \sy s, but it is worthwhile to recall that
they are so termed after Sinai, who in 1970 proved that they possess
the $K$-property. The same result was extended to higher dimension
by Sinai and Chernov in 1987 \cite{sc}. (A rather complete reference
on \bi s is \cite{t}. The hasty reader may find a brief history of
\emph{dispersing} \bi s in the introduction of \cite{l1}.)

\skippar

Turning to LGs, it is not surprising that much less work has been done
on them: the relevant invariant \me\ for these \sy s is infinite, and
ordinary \erg\ theory does not apply. Still it makes sense to pose the
question of the stochastic properties; only, one must utilize a more
refined \erg\ theory (as found, e.g., in the beautiful book by
Aaronson \cite{a}). For instance, the most sound definition of \erg
ity in infinite \me\ is the \emph{indecomposability} of the \sy\ 
modulo the \me. That is, a dynamical \sy\ endowed with an invariant
\me\ $\mu$ is \erg\ when the only invariant subsets mod $\mu$ have
\me\ zero or are complements of zero-\me\ sets (cf.\ Definition
\ref{def-erg} later on).

Another question that arises immediately for infinite-\me\ \sy s is
that of recurrence (cf.\ Definition \ref{def-rec}), because the
Poincar\'e Recurrence Theorem does not hold. For the PLG this has
baffled mathematicians for decades, and was answered in the
affirmative only in the last few years, by Schmidt \cite{sch} and
Conze \cite{co} independently. Both proofs use the Central Limit
Theorem (for a suitable class of observables) obtained by Bunimovich
and Sinai in 1981 \cite{bs}. All these results assume \emph{finite
horizon}: that is, the free path between collisions must be bounded
above.

Recurrence implies \erg ity (in the above sense) by a result of
Sim\'anyi \cite{si}.

We finish this short review of the PLG by mentioning a few other
important properties that are not so central in the present work.
Again for finite-horizon gases, the dynamics is diffusive and
converges to a Brownian motion upon suitable rescaling \cite{bs} (see
also \cite{bsc}.) The finite-horizon condition appears to be
physically relevant, as infinite-horizon PLGs are believed to have a
super-diffusive behavior (the mean square displacement at time $t$
increasing like $t \log t$) \cite{b}. A promising generalization of
the Central Limit Theorem is contained in a recent work by Sz\'asz and
Varj\'u \cite{sv}.

\skippar

In this paper we treat the generic LG with finite horizon. Subject to
very mild conditions on the curvature of the \sca s and their relative
position (cf.\ (\ref{cond-k})-(\ref{cond-tau})), we prove that
recurrence implies \erg ity (Theorem \ref{thm-erg}) and we give
several examples of recurrent aperiodic LGs (Propositions
\ref{prop-fin-mod} and \ref{prop-inf-mod}).

The exposition is thus organized: In Section \ref{sec-prel} we lay
down the necessary definitions and basic facts about \bi\ dynamics.
In Section \ref{sec-hyp} we show that our \sy s are \hyp.  In Section
\ref{sec-erg} we give the \erg ity result. In Section \ref{sec-rec} we
prove recurrence for a class of \sy s.

Throughout the paper we use the theory of \hyp\ \sy s with
singularities without really reviewing it for the reader. Rather, we
try to be as precise as possible as to what part of the theory we are
applying, or we need to modify, at any given time. Standard references
include \cite{cfs, sc, lw}.

\bigskip

\noindent
\textbf{Acknowledgments.} My thanks to N.~Chernov, N.~Sim\'anyi and
D.~Sz\'asz for having shared with me their expertise on this subject.
This work was partially supported by COFIN--MIUR, Italy (project: {\em
``Sistemi dinamici classici, quantistici e stocastici''}).

\sect{Preliminaries}
\label{sec-prel}

We study the \bi\ dynamics on $\R^2 \setminus \cup_{\a \in \is}
\Sc_\a$ (notation introduced in Section \ref{sec-intro}): A material
point moves with constant velocity until it hits the boundary of its
allowed domain. The point is then instantaneously scattered by the
obstacle according to the Fresnel law: the angle of incidence equals
the angle of reflection. It is clear that in this process the modulus
of the velocity never changes (a consequence of the conservation of
energy for this singular Hamiltonian \sy). So this modulus is
conventionally fixed to be 1. The dynamical \sy\ thus defined is
called the \emph{\bi\ flow} on $\R^2 \setminus \cup_{\a \in \is}
\Sc_\a$.

Since for most of the time the dynamics is trivial, we take the point
of view of discrete dynamical \sy s, and consider the Poincar\'e
cross-section corresponding to collisions with the obstacles. More
precisely, we restrict our phase space to points of the type $\x =
(q,v)$, with $q \in \partial \Sc_\a$, for some $\alpha$, and $v \in
\ts_q \R^2 = \R^2$, such that $|v| = 1$ and $v$ points outwardly with
respect to $\Sc_\a$. Any such point is also called a \emph{line
  element} and the whole set is denoted by $\ps$. Given an $\x \in
\ps$ we define $T \x =: \x_1 =: (q_1,v_1)$ to be first line element,
in the forward (flow-)\tr y of $\x$, that belongs to $\ps$ (see
Fig.~\ref{faplg1}). $T$ is then the Poincar\'e map induced by the flow
on the cross-section $\ps$, and the $T$-\o\ of $\x$ is the sequence of
snapshots of the flow-\tr y of $\x$, taken immediately after each
collision. (In this work I will try to be as consistent as possible
and use the term `\o' for an \o\ of a map, such as $T$; and `\tr y'
for an \o\ of the \bi\ flow.)

\fig{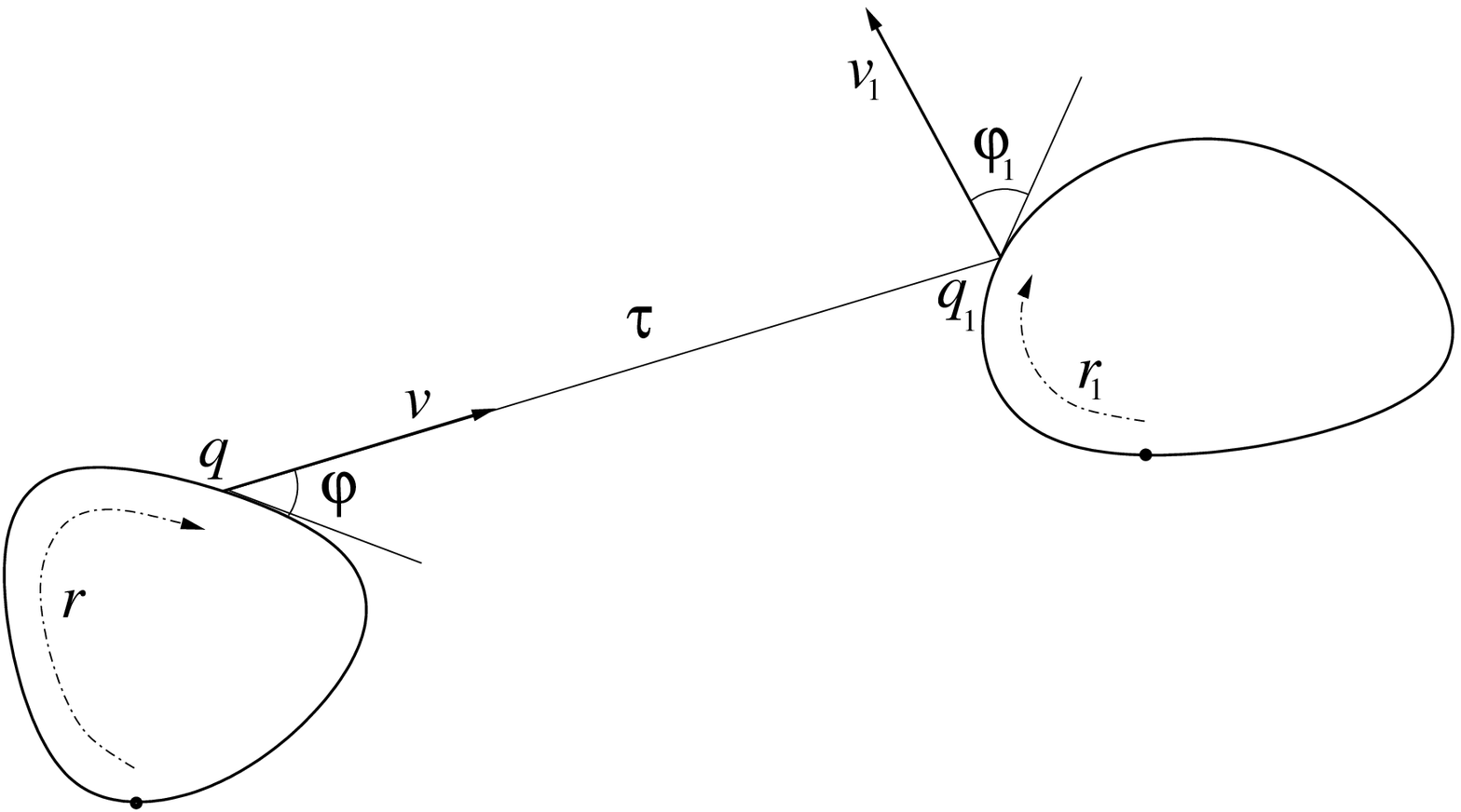} {4.7in} {Basic definitions for the \bi\ map.}

A convenient way to parametrize the phase space is to consider, for
each $\a$, the set $\ps_\a := S_{L_\a} \times [0,\pi]$, where $S_L$ is
the circle of circumference $L$ (i.e., the interval $[0,L]$ after
identification of the endpoints). A pair $(r,\ph) \in \ps_\a$
represents the unit vector, based on $q \in \partial \Sc_\a$, that
forms an angle $\ph$ with the clockwise tangent vector to $\partial
\Sc_\a$ in $q$; and $q$ is determined as the point whose clockwise
distance along $\partial \Sc_\a$ (from a certain fixed origin) is
$r$---see again Fig.~\ref{faplg1}. $L_\a$ is thus the length of the
closed curve $\partial \Sc_\a$. With a forgivable abuse of notation,
then, the phase space is rewritten as $\ps = \sqcup_{\a \in \is}
\ps_\a$ (the symbol denoting disjoint union).

We endow $\ps$ with the \me\ $\mu$ defined by $d\mu(r,\ph) = \sin\ph
\, dr d\ph$. This is the physically relevant \me\ because it is
induced on $\ps$ by the Liouville \me\ for the \bi\ flow (which, by
the way, is the uniform \me\ both in the position and in the momentum
coordinates). Hence $\mu$ is preseved by $T$ \cite{cfs}. Obviously
$\mu(\ps) = \infty$. From a mathematical point of view, this is
perhaps the most significant feature of the Lorentz gas.

It is well known that \sy s of the \bi\ type are discontinuous: any
line element $\x$ whose \tr y is tangent to the next obstacle, say
$\Sc_\b$, is a point of discontinuity for $T$, because other line
elements arbitrarily close to $\x$ may hit or miss $\Sc_\b$, which
causes $T \x$ to end up in completely different regions of the phase
space. The discontinuity set is then given by $\si := \si^{+} :=
T^{-1} \partial \ps$. Analogously, $\si^{-} := T \partial \ps$ is the
discontinuity set of $T^{-1}$. It is not hard to deduce that the
differential of $T$ blows up at $\si$ (see below). Therefore $\si$ is
usually called the \emph{singularity set}. It is made up of smooth
curves (\emph{singularity lines}) with the property that the endpoints
of any curve belong to other such curves \cite{s,sc,lw}.

For all elements in $\ps \setminus \si$, the differential of $T$ is
known (\cite[\S 14]{lw}, \cite[\S 3]{l1}) and, however, easily
computed with the help of Fig.~\ref{faplg1} and a little
patience. Denote by $\tau = \tau(\x) = \tau(r,\ph)$ the \emph{free
path} of the line element $\x$. Also, if $\x \in \ps_\a$, call
$k_\a(r)$ the curvature of $\partial \Sc_\a$ in the base-point of
$\x$: this is a smooth \fn\ by hypothesis, and positive by
convention. Also set $(r_1, \ph_1) := \x_1 = T \x$; $k = k_\a (r)$;
$k_1 = k_\b (r_1)$ (if $\x_1 \in \ps_\b$). Then
\begin{equation}
        D T_\x = \left[ 
        \begin{array}{cc} 
                \ds -\frac{\sin \ph} {\sin \ph_1} - \frac{k \tau}  
                {\sin \ph_1} & 
                \ds \frac{\tau} {\sin \ph_1} \ \ \\ \\
                \ds k + k_1 \frac{\sin \ph} {\sin \ph_1} + 
                \frac{ k k_1 \tau} {\sin \ph_1} \ \ &
                \ds -1 - \frac{k_1 \tau} {\sin \ph_1}
        \end{array}
        \right]. 
        \label{dmap}
\end{equation}

Before introducing the class of \sy s that we treat in this work
we need a very important, although standard, definition.

\begin{definition}
        The dynamical \sy\ $(\ps, T, \mu)$ is said to be
        \emph{recurrent} (in the sense of Poincar\'e) if, for every
        measurable set $A$, the \o\ of $\mu$-almost every $\x \in A$
        returns to $A$ at least once (and thus infinitely many times,
        due to the invariance of $\mu$). 
        \label{def-rec}
\end{definition}

We need not require the above property to hold for both the forward
and backward \o s because our \sy s are reversible; i.e., there exists
a map $I$ on $\ps$ such that $I^2 = id$ and $T^{-1} I T = I$. (This
map is of course $I(r,\ph) = (r,\pi-\ph)$.)

Since the Poincar\'e Recurrence Theorem fails here ($\mu(\ps) =
\infty$), one cannot take Definition \ref{def-rec} for granted and
indeed there are obvious cases of non-recurrent \sy s---think of a LG
with finitely many \sca s. (In fact these open \bi s are usually
studied only on the recurrent part of their phase space \cite{lm,
st}.)

\skippar

Let us accept the abuse of notation whereby $k(\x)$ means $k_\a(r)$,
for $\x = (r,\ph) \in \ps_\a$. Then we denote by $\mathcal{X}$ the
class of all \emph{recurrent} LGs for which there exist four positive
numbers, $k_m, k_M, \tau_m, \tau_M$, such that, $\forall \x \in \ps$,
\begin{eqnarray}
        && k_m \le k(\x) \le k_M;               \label{cond-k} \\ 
        && \tau_m \le \tau(\x) \le \tau_M.      \label{cond-tau}
\end{eqnarray}
The first condition is of a physical nature: we want to model a gas,
something that does not look much different in different parts of the
plane. In fact, this condition and the convexity of the \sca s imply
that no $\Sc_\a$ can be too big or too small; more precisely 
\begin{eqnarray}
        & 2 k_M^{-1} \le diam(\Sc_\a) \le 2 k_m^{-1}; &
        \label{diam-sca} \\
        & 2 \pi k_M^{-1} \le L_\a \le 2 \pi k_m^{-1}. &
        \label{len-sca}
\end{eqnarray}

Condition (\ref{cond-tau}) is primarily dictated by mathematical
feasibility, as we will see later. However, the existence of an upper
bound for the free path (which goes by the name of
\emph{finite-horizon condition}) does make a difference in the
dynamics, as recalled in the Introduction.

\skippar

For a given $\x$, define $\co(\x) := \rset {(dr, d\ph) \in \ts_\x \ps}
{dr \, d\ph \le 0}$, that is, the second and fourth quadrant of $\ts_\x
\ps$ in the $\{ \partial / \partial r, \partial / \partial \ph \}$
basis. This is usually called the \emph{unstable cone} at $\x$
\cite{lw, l1}. Take $u = (dr, d\ph) \in \co(\x)$. Any such vector is
called \emph{unstable}. The checkered sign configuration of
(\ref{dmap}) implies that $DT_\x u =: (dr_1, d\ph_1) \in \co(T \x)$,
and
\begin{eqnarray}
        (dr_1)^2 &\ge& \left( \frac{\sin \ph} {\sin \ph_1} \right)^2
        \left( 1 + \frac{k \tau} {\sin \ph} \right)^2 dr^2; \\ 
        (d\ph_1)^2 &\ge& \left( 1 + \frac{k_1 \tau} {\sin \ph_1}
        \right)^2 d\ph^2. 
\end{eqnarray}
Therefore, if we define the metric
\begin{equation}
        \| u \|_{\x}^2 = \| (dr, d\ph) \|_{(r,\ph)}^2 := \sin^2 \ph \:
        dr^2 + d\ph^2, 
        \label{inc-norm}
\end{equation}
and set $\lambda := 1 + k_m \tau_m > 1$, we obtain from
(\ref{cond-k})-(\ref{cond-tau}):
\begin{equation}
        \| DTu \|_{\x_1} \ge \lambda \, \| u \|_{\x}.
        \label{inc-lambda}
\end{equation}
We name (\ref{inc-norm}) the \emph{increasing norm} for unstable
vectors.

\sect{Hyperbolicity}
\label{sec-hyp}

In this section we see how to construct local stable and unstable
manifolds (LSUMs) for a.e.\ point of $\ps$, and how to prove that they
are absolutely continuous w.r.t.\ $\mu$. Although the subject is very
standard, it is useful for later purposes to state the definition of
LSUM. 

For $n \in \N \cup \{ \infty \}$, let us denote $\si_n^{\pm} :=
\bigcup_{i=0}^{n-1} T^{\mp i} \si^{\pm}$.

\begin{definition}
        Given a point $\x \in \ps \setminus \si_\infty^{+(-)}$, we
        define a {\em local (un)stable manifold $\wsu$ for $T$ at
        $\x$} to be a $C^1$ curve containing $\x$ in its
        (one-dimensional) interior, not intersecting
        $\si_\infty^{+(-)}$, and such that:
        \begin{itemize} 
                \item[(a)] The tangent spaces to $\wsu$ are included 
                in the pull-backwards (or push-forwards) of all future 
                (past) stable (unstable) cones. (E.g., in the unstable 
                case, this means: $\forall \y \in \wu$, $\ts_\y \wu 
                \subset \bigcap_{n\le 0} DT_{T^n \y}^{-n} \co(\y)$);

                \item[(b)] $\forall \y \in \wsu$, $d(T^n \y, T^n \x)
                \to 0$, as $n \to +\infty(-\infty)$;

                \item[(c)] If $\wsu_0$ is another such manifold, then
                so is $\wsu \cap \wsu_0$. 
        \end{itemize}
        \label{def-lsum}
\end{definition}

Condition \emph{(c)} is a sort of uniqueness property. For this reason
we refer to any LSUM as \emph{the} LSUM at $\x$ and denote it by
$\wsu(\x)$.  One readily sees that $T^{-1} \wu(\x)$ is a LUM at
$T^{-1} \x$ and $T \ws(\x)$ is a LSM at $T\x$. This earns the two
collections of local stable and unstable manifolds the name of
\emph{local invariant foliations}.

Condition $\emph{(a)}$ states that the tangent vectors to, say,
$\wu(\x)$ are unstable `at all orders'. This and (\ref{inc-lambda})
imply that if $\y \in \wu(\x)$ then, $\forall n \le 0$, $d_\|(T^n \y,
T^n \x) \le d_\|(\y, \x) \lambda^{n}$, where $d_\|$ is the distance in
the $\| \,\cdot\, \|$-norm.  We therefore say that $\wu(\x)$ is
\emph{exponentially unstable} in the increasing metric. We will see
later that this property holds for the Riemannian metric too.

\skippar

We do not recall the other basic notions of hyperbolic theory, such as
the definition of absolutely continuous foliation. (The reader who
would like to review the fundamentals of Pesin's theory, especially in
the context of \bi s, can consult \cite{ks,sc} and the references
therein.)

Even as concerns proofs, we only show the details of a couple of
lemmas that are specific to our \sy, and simply sketch the remaining
customary arguments---mainly to convince the reader that the
infiniteness of $\mu$ is no big trouble. It will also be noticed that
the recurrence property is \emph{not} needed here.

\skippar

The first observation one has to make to prove that a certain \bi\ is
\hyp\ is that the singularities are not a big obstruction. The lemma
that does the job is the following.

\begin{lemma}
        For a.e.\ $\x \in \ps$, these exists a $C_0 = C_0(\x)$ such that,
        $\forall n \ne 0$, 
        \begin{displaymath}
                d_\| (T^n \x, \si^{+} \cup \si^{-} \cup \partial \ps) 
                \ge C_0 \, |n|^{-4}.
        \end{displaymath}
        \label{lemma-appr}
\end{lemma}

This is in turn based on the next lemma, whose formulation needs two
extra definitions. Let $\si_\a^{\pm} := \si^{\pm} \cap \ps_\a$. Also,
for $A \subseteq \ps$ and $\eps>0$, set
\begin{equation}
        A_{(\eps)} := \rset{\x\in\ps} {d_\|(\x,A) \le \eps}.
        \label{a-eps}
\end{equation} 

\begin{lemma}
        There exists a constant $C_1$, \emph{independent of $\a \in 
        \is$}, such that, for $\eps>0$ small enough,
        \begin{displaymath}
                \mu \left( ( \si_\a^{+} \cup \si_\a^{-} \cup 
                \partial \ps_\a )_{(\eps)} \right) \le C_1 \, \eps.
        \end{displaymath}
        \label{lemma-sing}
\end{lemma}

\proofof{Lemma \ref{lemma-sing}} Obviously we can treat
$(\si_\a^{\pm})_{(\eps)}$ and $(\ps_\a)_{(\eps)}$ separately. For the
latter case the assertion is immediate since $(\ps_\a)_{(\eps)} = \{
(r,\ph) \in \ps_\a \,|\, \ph \le \eps \mbox{ or } \ph \ge \pi - \eps
\}$ and $L_\a$ is bounded above. (Remember that $\ps_\a$ is a cylinder
based on the circle of length $L_\a$.)

As for $\si_\a^{+}$, it is easy to see that it is composed of a
finite number of strictly increasing smooth curves $\gamma$, and this
number does not depend on $\a$. In fact, any \sca\ ``seen'' from
$\Sc_\a$ can only generate two tangent elements $(r,\ph)$ for any
given $r$, and, since $\Sc_\a$ is convex, $\ph$ must grow as $r$
increases (cf.\ Fig.~\ref{faplg1}). But the number of accessible \sca
s (we will also call them \emph{nearest neighbors}) is bounded above,
due to the finite-horizon condition and (\ref{diam-sca}). Notice that
the length of each $\gamma$ cannot exceed $L_\a + \pi$, which is less
than some universal constant by (\ref{len-sca}).

Analogous reasoning works for $\si_\a^{-}$, except that the curves are
strictly decreasing. 

So we only need prove that $\mu (\gamma_{(\eps)}) \le C_2 \, \eps$,
with $C_2$ independent of $\a$. This is essentially contained in
\cite[Lemma 2]{sc} but since the terminology used there is rather
different from that of the present work, we give a simple
``analytical'' proof here.

First of all, if $\gamma$ keeps away from $\partial \ps_\a$, the
result is obvious as $d_\|$ is equivalent to the Riemannian distance
$d$ there. So it suffices to look at, say, $\gamma \cap \{ \ph \le
\delta \}$, for some $\delta>0$. Since $k_\a (r)$ is bounded below,
$\gamma$ can only intersect $\{ \ph = 0 \}$ transversally. Therefore,
without loss of generality, we can focus on the case $\gamma = \{ 0 \}
\times [0, \delta]$. Given $\x = (r,\ph) \in \ps_\a$, the ball in the
unstable metric $B_\| (\x,\eps)$ is more or less an ellipse of
semiaxes $\eps / \sin\ph$ and $\eps$; and is certainly contained in
the rectangle $[ r - \eps / \sin (\ph-\eps) , r + \eps / \sin
(\ph-\eps) ] \times [ \ph - \eps, \ph + \eps ]$ (notice that this
might exceed $\ps_\a$).

\fig{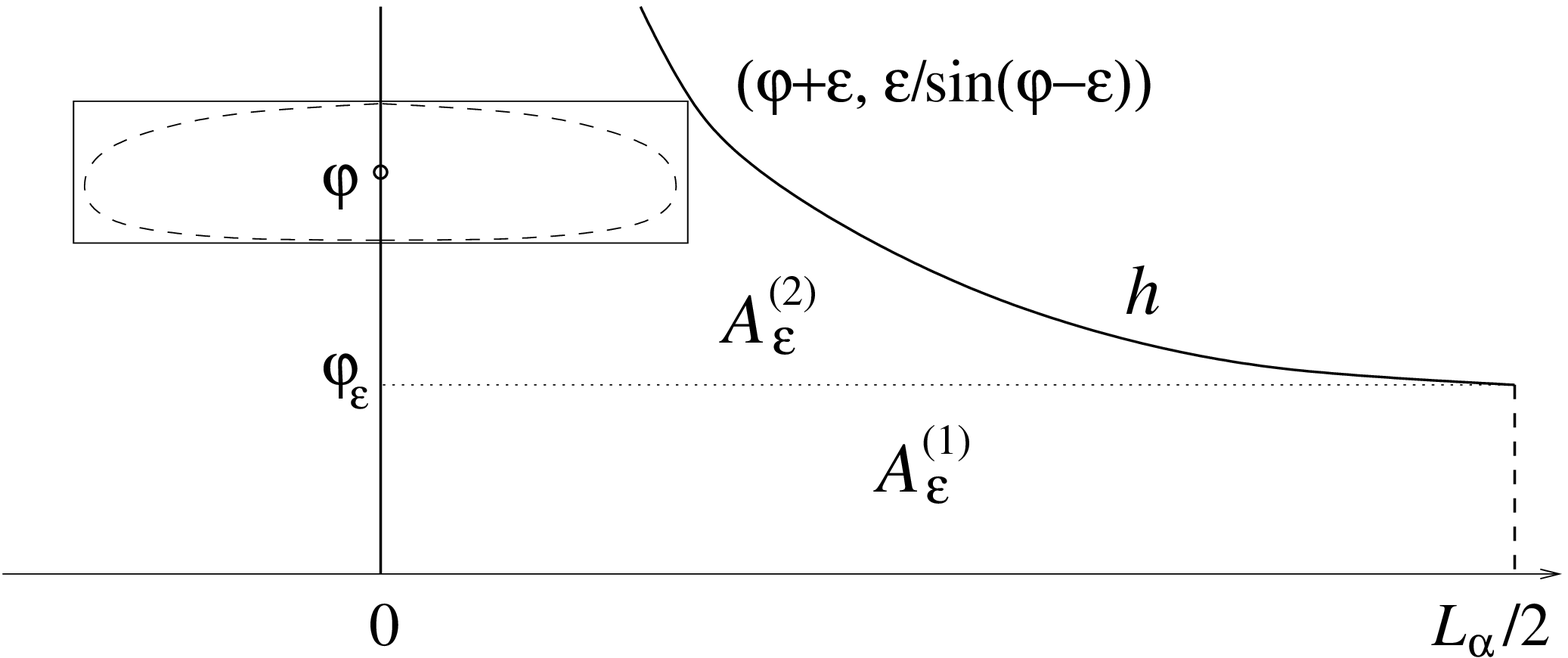} {4.2in} {The set $A_\eps$ as in the proof of Lemma
\ref{lemma-sing}.}

Hence, with the help of Fig.~\ref{faplg2}, we see that the ``right''
part of $\gamma_{(\eps)}$ is contained in $A_\eps := \rset{(r,\ph) \in
\ps_\a} {0 \le \ph \le \delta; \ 0 \le r \le h(\ph)}$, with
\begin{equation}
        h(\ph) = \left\{ 
        \begin{array}{cr}
                \ds \frac{\eps} {\sin (\ph - 2\eps) } & \quad \ph \in 
                [2\eps, \delta]; \\ 
                \infty & \quad \ph \in [0, 2\eps).
        \end{array}
        \right.
\end{equation}
But this is evidently an overestimate, as the right part of
$\gamma_{(\eps)}$ only comprises values of $r$ within $[0,
L_\a/2]$. Set $\ph_\eps$ so that $h( \ph_\eps ) = L_\a / 2$. This
means that
\begin{equation}
        \ph_\eps = \arcsin \left( \frac{2\eps} {L_\a} \right) + 2\eps
        \simeq 2\eps \left( \frac1 {L_\a} + 1 \right).
        \label{sing-20}
\end{equation}
(We say that $a \simeq b$ if both quantities depend on $\eps$ and $a/b
\to 1$, as $\eps \to 0^{+}$.) Discarding the exceeding part of
$A_\eps$, we split the reminaing part in two pieces:
\begin{eqnarray}
        A_\eps^{(1)} &=& [0, L_\a/2] \times [0, \ph_\eps]; \\
        A_\eps^{(2)} &=& \{ \ph_\eps \le \ph \le \delta; \ 0 \le r \le 
        h(\ph) \}. 
\end{eqnarray}
One can estimate $\mu( A_\eps^{(1)} )$ from above using the Lebesgue
\me; (\ref{sing-20}) and (\ref{len-sca}) ensure that this is less than
some $C_3 \, \eps$. On the other hand, for $\delta$ small,
\begin{eqnarray}
        \frac{ \mu (A_\eps^{(2)}) } {\eps} &=& \int_{\ph_\eps}^\delta
        \frac{\sin \ph} {\sin (\ph - 2\eps)} \, d\ph \simeq
        \int_{\ph_\eps}^\delta \frac{\ph} {\ph - 2\eps} \, d\ph =
        \int_{\ph_\eps - 2\eps}^{\delta - 2\eps} \left( 1 +
        \frac{2\eps} y \right) dy = \nonumber \\
        &=& \delta - \ph_\eps + 2\eps [ \log(\delta - 2\eps) - 
        \log( \ph_\eps - 2\eps ) ] \le 2\delta, 
\end{eqnarray}
as $\eps \to 0^{+}$, since from (\ref{sing-20}) $\ph_\eps - 2\eps
\simeq 2\eps / L_\a$.
\qed

\proofof{Lemma \ref{lemma-appr}} Without loss of generality we will
only consider forward semi\o s. Denote
\begin{equation}
        \nh_n(\a) := \rset{\b \in \is} { \exists \x \in \ps_\a, \:
        k \le n, \mbox{ such that } T^k \x \in \ps_\b }.
        \label{appr-10}
\end{equation}
In other words, $\nh_n(\a)$ is the index set of all
\emph{$k^{th}$-nearest neighbors} of $\Sc_\a$, for $k \le n$. By the
same line of reasoning as in the previous proof, it is clear that 
\begin{equation}
        \# \nh_n(\a) \le C_4 \, n^2,
        \label{appr-20}
\end{equation}
$C_4$ not depending on $\a$. Consider an $\x \in \ps_\a$. If the
inequality $d_\| (T^n \x, \si^{+} \cup \si^{-} \cup \partial \ps) \le
n^{-4}$ is verified only for a finite number of $n$'s, then $C_0(\x)$
can be found so that $d_\| (T^n \x, \si^{+} \cup \si^{-} \cup \partial
\ps) \ge C_0 \, n^{-4}$ foa all $n>0$. This does \emph{not} happen if,
and only if, $\x \in \ps_\a$ and
\begin{equation}
        T^n \x \in \bigcap_{m\ge 1} \: \bigcup_{n\ge m} 
        \left( \left( \si^{+} \cup \si^{-} \cup \partial \ps
        \right)_{(n^{-4})} \right). 
        \label{appr-30}
\end{equation}
But the arguments above about the phase-space accessibility of \o s
starting in $\ps_\a$ show that (\ref{appr-30}) is equivalent to
\begin{equation}
        \x \in \ps_\a \cap \, \bigcap_{m\ge 1} \: \bigcup_{n\ge m}
        T^{-n} \left( \bigcup_{\b \in \nh_n(\a)} \left( \si_\b^{+}
        \cup \si_\b^{-} \cup \partial \ps_\b \right)_{(n^{-4})}
        \right). 
        \label{appr-40}
\end{equation}
The \me\ of the above r.h.s.\ is zero by Borel-Cantelli, as the \me\
of the individual terms decreases like $n^{-2}$ (by Lemma
\ref{lemma-sing}, (\ref{appr-20}), and the invariance of $\mu$ w.r.t.\
$T$). Taking the union over $\a \in \is$ completes the proof.  
\qed

The existence and absolute continuity of the local invariant
foliations is a local matter. More precisely, for $\x \in \ps$, it has
to do with the behavior of $DT$ on a sequence of balls centered at
$T^n \x$ and of shrinking radii, as $|n| \to \infty$.

Say we want to consider LUMs. Lemma \ref{lemma-appr} implies that, for
some $C_5>0$, no element of the sequence $\{ B_\| (T^n \x, C_5
\lambda^{-|n|}) \}_{n\le 0}$ can intersect the singularities or the
boundary of the phase space. So, inside those balls, we have a
\emph{bona fide} Pesin's theory. The invariance of the unstable cones
($DT \co(\x) \subset \co(T\x)$) and (\ref{inc-lambda}) provide the
streching mechanism for \emph{unstable curves} (curves whose tangent
vectors are everywhere unstable). And the fact that $DT^{-n} \co(T^n
\x)$ collapses to a line (in this case, an easy consequence of
(\ref{inc-lambda}) and the invariance of $\mu$) proves that the
process of pushing forward unstable curves ``from the past'' has a
limit, which is precisely the LUM $\wu(\x)$.  To make things even
simpler, we are dealing with a \emph{uniformly hyperbolic} \sy\
(w.r.t.\ $\| \,\cdot\, \|$) and standard results \cite[\S6]{kh} show
that $\wu(\x)$ is as regular as $T$ is (meaning, away from $\si^{+}$).

Condensing in the same way the arguments for the absolute continuity
is impossible. Suffice it to say that the major task is to control the
\emph{distorsion coefficient} of $T$ along the sequence of shrinking
balls. The distorsion coefficient is a certain logarithmic derivative
of $|DT u|$, for some vector $u$, with $|u| = 1$; and the requirement
is that it grows less than exponentially as $n \to -\infty$. General
results on finite-horizon \bi s show that this quantity can increase
at most like a negative power of $\sin \ph_n$, with $(r_n, \ph_n) :=
\y_n \in B_\| (T^n \x, C_5 \lambda^n)$ (we are simplifying a bit
here---see \cite[\S7]{l1}).  But $\sin\ph_n \ge C_6 \min \{ \ph_n,
\pi-\ph_n \} = C_6 \, d_\|( \y_n, \partial \ps)$ and Lemma
\ref{lemma-appr} says in particular that this distance does not
approach zero faster than $|n|^{-4}$. Which verifies the requirement.

\begin{remark}
        The above also proves that the LSUMs are exponentially
        (un)stable w.r.t.\ the Riemannian metric. In fact $| \,\cdot\,
        |_{(r,\ph)} \le (\sin \ph)^{-1} \, \| \,\cdot\, \|_{(r,\ph)}$
        and in Definition \ref{def-lsum}, \emph{(b)} the polynomial
        growth of $(\sin \ph)^{-1}$ is tamed by the exponential
        contraction of $d_\|(T^n \y, T^n \x)$. Hence the convergence
        rate is any $\lambda' < \lambda$.
\end{remark}

None of the above assertions depend on the infiniteness of $\mu$, or
the non-compact\-ness of the \bi. To describe it in lay terms, the
particle has no way of knowing whether it is traveling in a plane or
in a torus, by only looking at small neighborhoods of its \o.

\sect{Ergodicity} 
\label{sec-erg}

For dynamical \sy s with an infinite invariant \me, the many available
definitions of \erg ity fail to be equivalent. The definition that is
usually retained is the following \cite{a}:

\begin{definition}
        A dynamical \sy\ $(\ps, T, \mu)$ is \erg\ if every $A
        \subseteq \ps$, measurable and invariant mod $\mu$ (i.e.,
        $\mu( T^{-1} A \, \Delta \, A) = 0$), has either zero \me\ 
        or full \me\ (i.e., $\mu (\ps \setminus A) = 0$).
        \label{def-erg}
\end{definition}

Definition \ref{def-erg} does not require the \sy\ to be recurrent,
but it usually makes little sense if there is a dissipative part (thus
causing the whole \sy\ to be dissipative). We do not have that problem
here.

The goal of this section is to prove the next result.

\begin{theorem}
        Let the Lorentz gas $\{ \Sc_\a \}_{\a \in \is}$ belong to
        $\mathcal{X}$. The associated dynamical \sy\ $(\ps, T, \mu)$
        is \erg.  
        \label{thm-erg}
\end{theorem}

The essential tool to work out Theorem \ref{thm-erg} is the \emph{local
\erg ity theorem}, also known as the \emph{fundamental theorem} for
\hyp\ \bi s. We state it in a rather technical guise that is
convenient for our purposes.

\begin{theorem}
        Let $A$ be a full-\me\ subset of $\ps$. Then, for any $\x_0$ that
        possesses a semi\o\ (i.e., $\x_0 \in \ps \setminus
        \si_\infty^{+}$, or $\x_0 \in \ps \setminus \si_\infty^{-}$),
        there is a neighborhood $U$ of $\x_0$ with the following
        property:

        A.e.\ two points $\x', \x'' \in U$ are connected by a finite
        alternating sequence of LSUMs, $\ws(\x_1), \wu(\x_2), ... ,
        \wu(\x_m)$, with $\x_1 := \x'$ and $\x_m := \x''$. Each LSUM
        intersects the next transversally in a point of $A$. Also
        $\x_1, ... , \x_m \in A$.
        \label{thm-loc-erg}
\end{theorem}

We do not give the complete proof of the local \erg ity theorem. As the
name suggests, most of it is strictly local and does not require $\mu$
to be finite \cite[\S\S8-12]{lw}. The technical conditions that are
called for are all verified for dispersing \bi s; this has been known
since \cite{s}.

The non-local part in the proof of Theorem \ref{thm-loc-erg} is in
general the most laborious and for this reason has its own name, the
\emph{tail bound} \cite[\S13]{lw}. We state it here in the form of a
lemma and prove it in detail. 

First we assume that the LSUMs are maximal; i.e., for every $\x$ for
which a LUM exists, we take $\wu(\x)$ to be the union of all the LUMs
at $\x$. This implies that, for a.e.\ $\x$, $\partial \wu(\x) \subset
\partial \ps \cup \si_\infty^{-}$.  (It might also happen that
$\partial \wu(\x)$ intersects the closure of $\si_\infty^{-}$---which
is typically the whole $\ps$---but standard reasonings show that this
cannot occur for more than a null-\me\ set of $\x$'s. See
\cite[Rk.~8.9]{l1}.)  The analogous statement holds of course for
$\ws(\x)$.

If $\y \in \partial \wu(\x) \cap T^m \si^{-}$, for some $m \ge 0$, we
say that $\wu(\x)$ is \emph{cut} by $T^m \si^{-}$ at $\y$.  The radius
of $\wu(\x)$ is defined as the minimum distance, \emph{along $\wu(\x)$
and in the increasing metric}, between $\x$ and the two endpoints of
$\wu(\x)$.

\begin{lemma}
        For every $\x_0 \in \ps$, there is a neighborhood $U_0$ of 
        $\x_0$, and a $\delta_0 > 0$ such that, $\forall \eta > 0$, 
        $\exists M$ that verifies
        \begin{displaymath}
                \mu \left( \rset{\x \in U_0} {\wu(\x) \mbox{ \emph{has
                radius}} < \delta \mbox{ \emph{because it is cut by}} 
                \bigcup_{m=M+1}^\infty T^m \si^{-} } \right) \le
                \eta \delta,
        \end{displaymath}
        for every $\delta \le \delta_0$.
        \label{lemma-tail}
\end{lemma}

\proofof{Lemma \ref{lemma-tail}} Let $Y = Y(U_0, \delta, M)$ be the set in
the statement of the lemma and, $\forall \x \in Y$, set $m(\x)$ to be
the smallest integer such $T^m \si^{-}$ cuts $\wu(\x)$ within a
distance $\delta$ of $\x$. If we define $Y_m := \rset{\x \in Y} {m(\x)
= m}$, then
\begin{equation}
        \mu(Y) = \sum_{m>M} \mu(Y_m) = \sum_{m>M} \mu(T^{-m} \, Y_m).
        \label{tail-10}
\end{equation}
The first equality holds because the $Y_m$'s are the level sets of the
\fn\ $m$; the second equality is due to the invariance of $\mu$.

For $\x \in Y_m$ denote by $\y = \y(\x) \in \partial \wu(\x) \cap T^m
\si^{-}$ the point where the cut occurs.  By the exponential
contraction of $\wu(\x)$, the portion of $T^{-m} \wu(\x)$ between
$T^{-m} \x$ and $T^{-m} \y \in \si^{-}$ has length less than $\delta
\lambda^{-m}$ (in the increasing metric). So $d_\|(T^{-m} \x, \si^{-})
< \delta \lambda^{-m}$.

Let $\a \in \is$ so that $\x_0 \in \ps_\a$, and select any $U_0
\subseteq \ps_\a$.  Following the line of reasoning of Lemma
\ref{lemma-appr}, $T^{-m} Y_m \subset \cup_{\b \in \nh_m(\a)} \,
\ps_\b$. This and the above imply that
\begin{equation}
        T^{-m} \x \in \bigcup_{\b \in \nh_m(\a)}
        (\si_\beta^{-})_{(\delta \lambda^{-m})}, 
        \label{tail-20}
\end{equation}
whose measure, by Lemma \ref{lemma-sing} and (\ref{appr-20}), is $\le
C_1 C_4 \, m^2 \, \delta \lambda^{-m}$. Therefore, from
(\ref{tail-10}),
\begin{equation}
        \frac{\mu(Y)} \delta \le C_1 C_4 \sum_{m>M}  m^2 \,
        \lambda^{-m} \to 0,
\end{equation}
as $M \to +\infty$.
\qed

\begin{lemma}
        Fixed an $\a \in \is$ and denoted by $T_\a$ the first return
        map induced by $T$ on $\ps_\a$ (this is well defined a.e.\ by
        recurrence), $(\ps_\a, T_\a, \mu)$ is \erg.  
        \label{lemma-erg}
\end{lemma}

\proof We proceed in three steps:

\skippar
\noindent
\textsc{Step 1:} {\itshape The LSUMs relative to $T$ are also LSUMs
relative to $T_\a$.} Let us see to it. In Definition \ref{def-lsum},
part \emph{(c)} has nothing to do with $T$ or $T_\a$. Part \emph{(a)}
is obvious since $T_\a^k \y = T^{n_k} \y$, for some $n_k < 0$. For
part \emph{(b)} it is sufficient to prove that if $\y \in \wsu(\x)$
then the (future or past) return times to $\ps_a$ are the same for
$\y$ and $\x$. But this follows from \emph{(b)} itself as, by
definition, $d_\| (\ps_\a, \ps_\b) = \infty$ when $\a \ne \b$;
therefore the sequence of \sca s hit is the same for for $\y$ and
$\x$.

\skippar
\noindent
\textsc{Step 2:} {\itshape Any $\x_0$ that possesses a semi\o\
has a neighborhood $U$ contained in one \erg\ component.}
This is Hopf's idea. Take a measurable $f$, continuous (and thus
uniformly continuous) on $\ps_\a$, and consider its forward and
backward \erg\ averages w.r.t.\ $T_\a$:
\begin{equation}
        f^{\pm}(\x) := \lim_{n \to +\infty} \frac1n \sum_{k=0}^{n-1}
        (f \circ T_\a^{\pm k})(\x).
\end{equation}
By a standard application of Birkhoff's Theorem, the set 
\begin{equation}
        A := \rset{\x \in \ps} {f^{+}(\x), f^{-}(\x) \mbox{ exist and
        coincide} }
\end{equation}
has full \me. Let us use it in Theorem \ref{thm-loc-erg}. In the
resulting neighborhood $U$, a.e.\ two points $\x', \x''$ are connected
by a curvilinear polyline made up of LSUMs. Call $\y$ the intersection
point between $\ws(\x_1)$ and $\wu(\x_2)$. From Step 1 and the uniform
continuity of $f$, $f^{+}(\y) = f^{+}(\x_1)$ and $f^{-}(\y) =
f^{-}(\x_2)$. But $f^{+}(\y) = f^{-}(\y)$ because $\y \in
A$. Proceeding this way we conclude that $f^{+}(\x') = f^{+}(\x'')$;
that is, $f^{+}$ is a.e.\ constant in $U$.  This continues to hold for
every $f \in L^1(\ps)$, via Birkhoff's Theorem and customary density
arguments.

\skippar
\noindent
\textsc{Step 3:} {\itshape There is only one \erg\ component for
$(\ps_\a, T_\a, \mu)$.} The set of points that fail to have a semi\o\
is $\si_\infty^{+} \cap \si_\infty^{-}$. But, as explained in Section
\ref{sec-hyp}, $\si_\infty^{+}$ and $\si_\infty^{-}$ are countable
unions of smooth curves, respectively strictly increasing and strictly
decreasing.  Therefore there can be at most one point of intersection
for each pair of increasing--decreasing curves; that is, at most
countably many points. Hence $\ps_\a \setminus ( \si_\infty^{+} \cap
\si_\infty^{-} )$ is path-connected.  
\qed

Now it is simple to prove the main result of Section \ref{sec-erg}.

\skippar

\proofof{Theorem \ref{thm-erg}} We claim that the \erg\ decomposition
of $T$ is coarser than the partition $\{ \ps_\a \}_{\a \in \is}$. In
fact the existence of $A, B$, invariant subsets, such that $\mu(A \cap
\ps_\a), \mu(B \cap \ps_\a) > 0$, for some $\a \in \is$, contradicts
Lemma \ref{lemma-erg}.

Therefore, if there were more than one \erg\ component, there would
be two nearest neighbors $\Sc_\a$ and $\Sc_\b$ such that $\ps_\a$ and
$\ps_\b$ are contained in two different invariant sets. But this is
absurd since one can always find $A \subset \ps_\a$, $\mu(A) > 0$,
such that $T(A) \subset \ps_\b$.
\qed

\sect{Recurrence}
\label{sec-rec}

We have recalled in the introduction that a finite-horizon periodic
Lorentz gas is recurrent (as in Definition \ref{def-rec})
\cite{co,sch} and thus belongs to $\mathcal{X}$ (conditions
(\ref{cond-k})-(\ref{cond-tau}) being trivially satisfied).
Here we present other examples from the class $\mathcal{X}$. 

\begin{definition}
        The LG $\{ \Sc_\a \}_{\a \in \is}$ is called a \emph{finite
        modification} of the finite-horizon PLG $\{ \Sc_\a \}_{\a \in
        \is_P}$ if $\is = (\is_P \setminus \is_1 ) \cup \is_2$, where:
        \begin{itemize} 
                \item[(a)] $\is_1$ is a finite subset of $\is_P$;

                \item[(b)] $\is_2$ is the index set of a finite LG
                such that $d (\Sc_\a, \Sc_\b) > 0$ for any $\a \in
                \is_2$, $\b \in \is_P \setminus \is_1$.
        \end{itemize}
        \label{def-fin-mod}
\end{definition}

\begin{remark}
        Calling the above a `perturbation' does not seem
        appropriate, as the modification can be arbitrarily big,
        provided it is finite.
\end{remark}

An example of a finite modification of a PLG is illustrated in
Fig.~\ref{faplg3}. Since it is clear that Definition \ref{def-fin-mod}
implies (\ref{cond-k})-(\ref{cond-tau}), we focus on recurrence.

\begin{proposition}
        A finite modification of a PLG is recurrent.
        \label{prop-fin-mod}
\end{proposition}

\proof Let us introduce $\ij := \rset{\b \in \is} {\b \in \nh_1(\a),
  \mbox{ for some } \a \in \is_2}$. This is $\is_2$ plus the indices
of the nearest neighbors to the \sca s labeled by $\is_2$ (cf.\ 
Fig.~\ref{faplg3}). Also denote $\ps_\ij := \sqcup_{\a \in \ij}
\ps_\a$. We proceed in three steps:

\fig{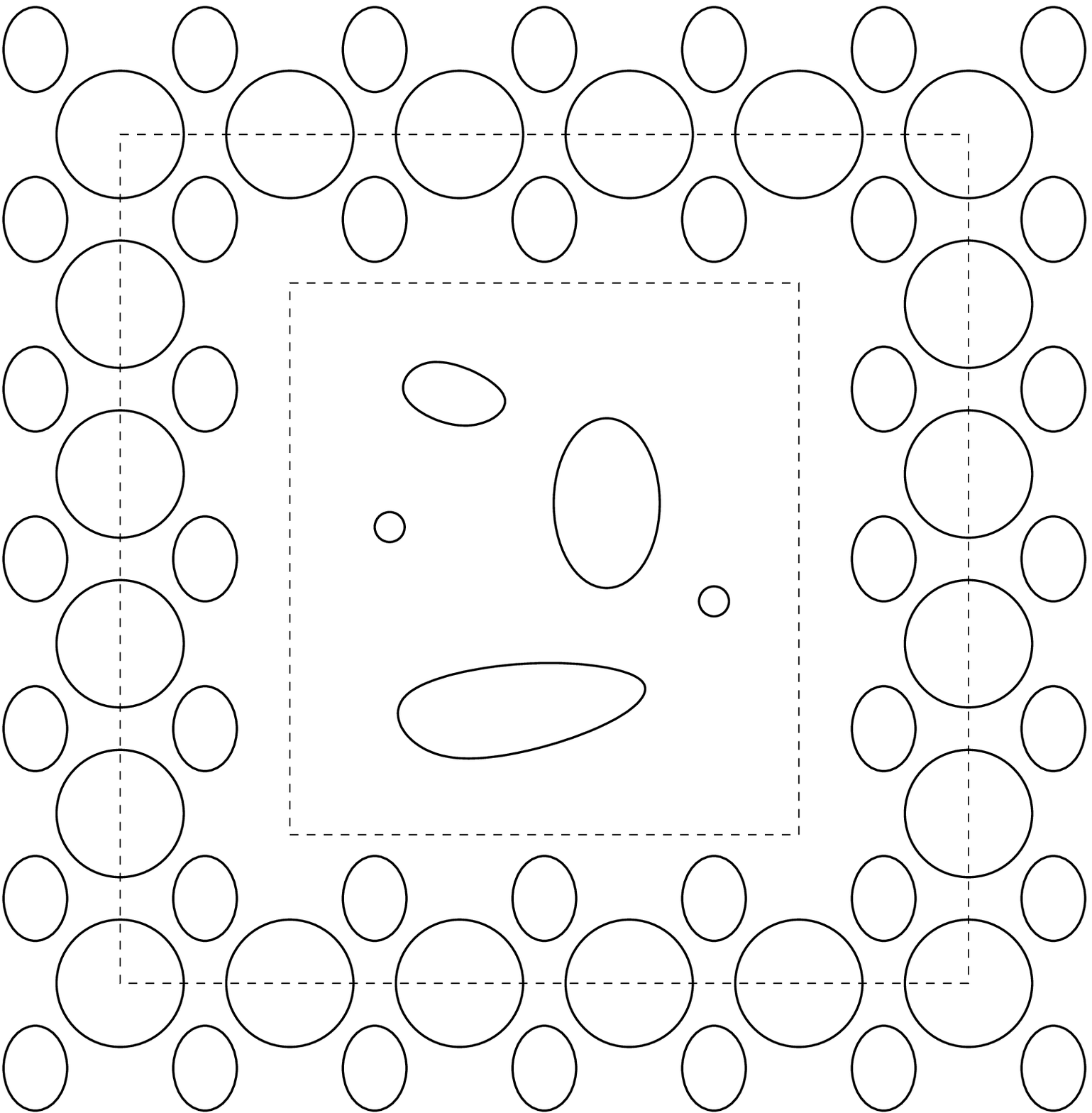} {4in} {A finite modification of a periodic Lorentz
gas. The inner dashed square encloses the \sca s associated to
$\is_2$; the outer square corresponds to $\ij$.}

\noindent
\textsc{Step 1:} {\itshape A.e.\ $\x \in \ps_\ij$ returns to
$\ps_\ij$.} First of all, for \emph{every} $\x \in \ps_{\is_2}$
(notation understood), $T\x \in \ps_\ij$. Then consider the set $\rset
{\x \in \ps_{\ij \setminus \is_2}} {T^n \x \not\in \ps_\ij, \, \forall
n>0}$. The forward \o\ of every point there lies the ``periodic part''
of the plane. One uses the recurrence of a PLG to conclude that this
set has \me\ zero.

\skippar
\noindent
\textsc{Step 2:} {\itshape Denoted by $T_\ij$ the return map to
$\ps_\ij$ (well defined a.e.\ by Step 1), $T_\ij$ is recurrent.}  This
is the Poincar\'e Theorem applied to $(\ps_\ij, T_\ij, \mu)$.

\skippar
\noindent
\textsc{Step 3:} {\itshape $T$ is recurrent.} It suffices to prove
that $\ps_\ij$ is a global cross-section (mod $\mu$) for $T: \ps
\longrightarrow \ps$. But a positive-\me\ subset of $\ps$ whose \o\
never intersects $\ps_\ij$ would contradict the \erg ity of a PLG.
\qed

If the \sy s above seem too close to those already known, we provide
further examples of LGs whose aperiodicity is not limited to a compact
region of $\R^2$.

\begin{proposition}
        There are LGs in the class $\mathcal{X}$ that are not finite
        modifications of PLGs.
        \label{prop-inf-mod}
\end{proposition}

\proof We present an algorithm to construct such gases
explicitly. Start with a PLG $\mathcal{L}_0 := \{ \Sc_\a \}_{\a \in
\is_0}$ and denote by $T_0$ the corresponding \bi\ map. Furthermore,
fix an index $\a \in \is_0$ and a point $O \in \Sc_\a$ which will be
our origin from now on.

For $\x \in \ps_\a$, denote by $n_1 = n_1(\x, T_0) \in \N \cup \{
+\infty \}$ the first return time of $\x$ to $\ps_\a$, w.r.t.\ $T_0$.
Indicated by $q(\x)$ the point in $\R^2$ corresponding to $\x$, we
define
\begin{equation}
        A_0(R) := \rset{\x \in \ps_\a} { \sup_{0<k<n_1} dist( q(T_0^k
        \x), O) \le R }, 
        \label{inf-mod-10}
\end{equation}
where $dist$ is the distance in $\R^2$. Since $T_0$ is recurrent,
$A_0(R) \nearrow \ps_a$ mod $\mu$, as $R \to +\infty$. Hence there is
an $R_0$ such that $\mu( \ps_\a \setminus A_0(R_0) ) \le \eps_0$. Here
$\{ \eps_k \}$ is a sequence of positive numbers such that $\eps_k
\searrow 0$, for $k \to +\infty$.

Now modify $\mathcal{L}_0$ in the annulus of radii $R_0$ and $R_0 +
\rho$ ($\rho$ is a sufficiently large fixed quantity). This means that
one can modify the shape, the position and the number of \sca s within
that annulus. Call the resulting LG $\mathcal{L}_1 := \{ \Sc_\a \}_{\a
\in \is_1}$, and the resulting map $T_1$. This is a finite
modification of a PLG. Furthermore, denoted by $T_{\a,k}$ the first
return map to $\ps_a$ relative to $T_k$, it is obvious that
$\left. \left( T_{\a,1} \right) \right|_{A_0(R_0)} = \left. \left(
T_{\a,0} \right) \right|_{A_0(R_0)}$.

It is clear that this process can be repeated recursively. Starting
with $\mathcal{L}_k$, a finite modification of a PLG, one defines
$A_k(R)$ as in (\ref{inf-mod-10}), with $T_k$ in the place of $T_0$.
By Proposition \ref{prop-fin-mod}, $A_k(R) \nearrow \ps_a$ mod $\mu$,
as $R \to +\infty$.  One then chooses $R_k > R_{k-1} + \rho$ such that
\begin{equation}
        \mu( \ps_\a \setminus A_k(R_k) ) \le \eps_k
        \label{inf-mod-20}
\end{equation}
and makes a modification in the annulus of radii $R_k$ and $R_k +
\rho$, to obtain $\mathcal{L}_{k+1}$ and $T_{k+1}$. We emphasize that
this modification can be made so that the curvature of the \sca s and
the free path among them stay bounded above and below, to comply with
(\ref{cond-k})-(\ref{cond-tau}).  Evidently, 
\begin{eqnarray}
        & A_k(R_k) = A_{k+1}(R_k) \subseteq A_{k+1}(R_{k+1}); & 
        \label{inf-mod-30} \\
        & \left. \left( T_{\a,k+1} \right) \right|_{A_k(R_k)}
        = \left. \left( T_{\a,k} \right) \right|_{A_k(R_k)}. &
        \label{inf-mod-40}
\end{eqnarray}
Since we never change the LG twice in the same region of the plane,
this modification process has a limit, which we call
$\bar{\mathcal{L}}$, with the corresponding maps, $\bar{T}$ for the LG
and $\bar{T}_\a$ for the returns to $\ps_a$.  In other words, from
(\ref{inf-mod-20})-(\ref{inf-mod-40}), $\bar{T}_\a: \ps_\a
\longrightarrow \ps_\a$ is well defined a.e.\ via $\left. \left(
\bar{T}_\a \right) \right|_{A_k(R_k)} = \left. \left( T_{\a,k} \right)
\right|_{A_k(R_k)}$.

By construction, \emph{all} points of $\cup_k A_k(R_k)$ come back to
$\ps_a$. So, as far as $\bar{T}$ is concerned, $\ps_a$ is contained in
$B$, the \emph{conservative} part of $\bar{\ps}$ (i.e., the complement
of its dissipative part \cite{a}). $B$ is obviously
$\bar{T}$-invariant.

We claim that the LSUMs of $\bar{T}$ that intersect $B$ are wholly
contained in $B$, mod $\mu$.  (The LSUMs exist a.e.\ and are
absolutely continuous by Section \ref{sec-hyp}---no recurrence was
used there.)  More precisely, this means that there is a set $B_0$,
$\mu( B \Delta B_0 ) = 0$, such that, called $B_1$ the union of all
the LSUMs based in $B_0$, we have $\mu( B \Delta B_1 ) = 0$.

In fact, take $\x \in B \cap \ps_\b$, $\b \in \bar{\is}$; passing
maybe to a full-\me\ subset $B_0$, we can assume that $\x$ returns to
$\ps_\b$ infinitely many times, both in the past and in the
future. Then any $\y \in \ws(\x)$, say, returns to $\ps_\b$ infinitely
many times in the future (Step 1 of Lemma \ref{lemma-erg}). Therefore,
save for a null-\me\ set of exceptions, $\y$ cannot be part of a
wandering set and thus belongs to $B$.

At this point, Theorem \ref{thm-loc-erg} proves that the decomposition
$\{ B, \bar{\ps} \setminus B \}$ is coarser than $\{ \ps_\a \}_{\a \in
\bar{\is}}$. That $B = \bar{\ps}$ mod $\mu$ follows from the
invariance of $B$ as in the proof of Theorem \ref{thm-erg}.  
\qed

\end{document}